%% file: graphicalBtienmt_R1_copy.tex
\documentclass[APA,Garamond1COL]{WileyNJDv5} 

\usepackage{amsmath,amssymb,bbm}

\newcommand{\eig}{\mathrm{eig}}

\newtheorem{assume}{Assumption}

\newtheorem{cor}{Corollary}[section]

\newtheorem{rmk}{Remark}
\DeclareMathOperator*{\argmin}{arg\,min}

\allowdisplaybreaks

\articletype{Article Type}%

\begin{document}

\title{Concentration of a sparse Bayesian model with Horseshoe prior in estimating high-dimensional precision matrix}

\author[1]{The Tien Mai}

\authormark{MAI \textsc{et al.}}
\titlemark{Concentration of a sparse Bayesian model with Horseshoe prior in estimating high-dimensional precision matrix}

\address{\orgdiv{Department of Mathematical Sciences}, \orgname{Norwegian University of Science and Technology}, \orgaddress{\state{Trondheim 7034}, \country{Norway}}}

\corres{\email{the.t.mai@ntnu.no}}

\abstract[Abstract]{Precision matrices are crucial in many fields such as social networks, neuroscience, and economics, representing the edge structure of Gaussian graphical models (GGMs), where a zero in an off-diagonal position of the precision matrix indicates conditional independence between nodes. In high-dimensional settings where the dimension of the precision matrix \( p \) exceeds the sample size \( n \) and the matrix is sparse, methods like graphical Lasso, graphical SCAD, and CLIME are popular for estimating GGMs. While frequentist methods are well-studied, Bayesian approaches for (unstructured) sparse precision matrices are less explored. The graphical horseshoe estimate by \cite{li2019graphical}, applying the global-local horseshoe prior, shows superior empirical performance, but theoretical work for sparse precision matrix estimations using shrinkage priors is limited. This paper addresses these gaps by providing concentration results for the tempered posterior with the fully specified horseshoe prior in high-dimensional settings. Moreover, we also provide novel theoretical results for model misspecification, offering a general oracle inequality for the posterior. A concise set of simulations is performed to validate our theoretical findings.
}



\maketitle

\section{Introduction}
\label{sc_intro}
Precision matrices play a pivotal role in numerous fields, including social networks, neuroscience, and economics \cite{pourahmadi2011covariance,ryali2012estimation,fan2016overview,callot2021nodewise}. The configuration of the precision matrix $\Omega = \Sigma^{-1} \in \mathbb{R}^{p \times p}$ represents the edge structure of a Gaussian graphical model, where the nodes are normally distributed as $\mathcal{N}(0, \Sigma)$ \cite{lauritzen1996}. Gaussian graphical models (GGMs) continue to be favored for network estimation due to the straightforward interpretation of the precision matrix: a zero in an off-diagonal position indicates conditional independence between the two corresponding nodes, given all other nodes. Due to these important implications, the estimation of precision matrices has received substantial attention over the past few decades.

In practical applications, users often encounter situations where the dimension $p$ of the precision matrix is comparable to or even exceeds the sample size $n$, and the precision matrix exhibits sparsity. In such contexts, leveraging the sparsity of the precision matrix is crucial to discern the edge structure of the graph and assess the conditional dependencies among the nodes. Among the most widely used frequentist methods for estimating GGMs are the graphical lasso \cite{friedman2008sparse}, the graphical SCAD \cite{fan2009network}, and the CLIME estimator \cite{cai2011constrained}. For a comprehensive overview, see \cite{cai2016estimating} and \cite{jankova2017honest}. These methods provide estimates of high-dimensional inverse covariance matrices under various sparsity patterns.
 
The body of work on Bayesian methodologies for \emph{unstructured precision matrices} is relatively limited. \cite{wang2012bayesian} introduced a Bayesian variant of the graphical lasso, utilizing block Gibbs sampling. \cite{banerjee2015bayesian} employed a similar prior structure to the Bayesian graphical lasso but included a substantial point-mass at zero for the off-diagonal elements of the precision matrix. They derived posterior convergence rates in the Frobenius norm under sparsity assumptions and proposed a Laplace approximation method for calculating marginal posterior probabilities of models. Additionally, \cite{gan2019bayesian} proposed spike-and-slab variants with double exponential priors. \cite{atchade2019quasi} combined the spike-and-slab prior with a pseudo-likelihood approach, resulting in a quasi-posterior distribution with favorable contraction properties.

A common issue with the spike-and-slab approach is the presence of binary indicator variables, which can hinder posterior exploration. The graphical horseshoe method, introduced by \cite{li2019graphical}, addresses these issues by applying the global-local horseshoe prior \cite{carvalho2010horseshoe} to Gaussian graphical models. \cite{li2019graphical} provided substantial empirical evidence of the graphical horseshoe estimate's superior performance compared to several competing Bayesian and frequentist methods in both low and high-dimensional settings. However, there is limited research on the theoretical properties of shrinkage priors for sparse precision matrix estimation. Recent studies by \cite{zhang2022contraction} and \cite{sagar2024precision} have adapted several versions of Horseshoe prior in different ways. Specifically, \cite{zhang2022contraction} consider a pseudo-likelihood approach similar to \cite{atchade2019quasi}, while \cite{sagar2024precision} explore a horseshoe-like prior. Moreover, both studies offer theoretical results only in an asymptotic sense and do not address high-dimensional settings. This work provides an attempt to address these gaps in the literature.

In this paper, we provide concentration results for the posterior in estimating precision matrices utilizing the horseshoe prior as in \cite{li2019graphical}. Following \cite{banerjee2015bayesian}, \cite{zhang2022contraction}, and \cite{sagar2024precision}, we assume that the underlying precision matrix is sparse. Importantly, our results are derived in a high-dimensional setting where \( p > n \). Moreover, while \cite{zhang2022contraction} and \cite{sagar2024precision} consider either a horseshoe-like prior or a horseshoe prior with a fixed global shrinkage parameter, our results are obtained for the fully specified horseshoe prior used in practice as in \cite{li2019graphical}. Another novel aspect of our work is that we also provide theoretical results for the case of model misspecification. Specifically, we prove a general oracle inequality for the posterior under an abstract setting that allows for model misspecification. A brief simulation study is conducted to validate our theoretical results.

The rest of the paper comprises the following sections
Section \ref{sc_modelmethod} present the problem of sparse precision matrix estimation, a Bayesian approach using tempered posterior and the Horseshoe prior. Section \ref{sc_concentration} contains our main theoretical results. Several simulations are conducted in Section \ref{sc_simu} to further illustrate and validate our theoretical findings. We conclude in Section \ref{sc_conclusion}. All technical proofs re gathered in Appendix \ref{sc_proofs}.

\section{Problem and method}
\label{sc_modelmethod}
Notations: 
Let $P,R$ be two probability measures and $\mu$ be any measure such that $P\ll \mu$ and $R\ll \mu$. For $\alpha\in(0,1)$, the $\alpha$-R\'enyi divergence  between  $P$ and $R$ is  defined by
$
D_{\alpha}(P,R)  =
\frac{1}{\alpha-1} \log \int \left(\frac{{\rm d}P}{{\rm d}\mu}\right)^\alpha \left(\frac{{\rm d}R}{{\rm d}\mu}\right)^{1-\alpha} {\rm d}\mu  \text{,}
$
and the Kullback-Leibler divergence is defined by
$
\mathcal{K}(P,R)  = 
\int \log \left(\frac{{\rm d}P}{{\rm d}R} \right){\rm d}P 
$
  if  $ P \ll R
$, and  $
+ \infty \text{ otherwise}.
$
 Let $ \|\cdot\|_q $ denote the $ \ell_q $-norm, $ \|\cdot\|_\infty $ denote the max-norm of vectors and let $ \|\cdot\|_0 $ denote the $ \ell_0 $ (quasi)-norm (the number nonzero entries) of vectors. .

\subsection{Model}

We are interested in estimating the precision matrix $\Omega \in \mathbb{R}^{p\times p}$ from the sample data $Y\in \mathbb{R}^{n\times p} $, where $n$ denotes the sample size, and $p$ is the dimension of precision matrix or equivalently the number of nodes in the corresponding graph. The rows of $Y$ are assumed to be i.i.d. $ \mathcal{N}( 0,\Omega^{-1})$ distributed. The likelihood function based on the data can be written as
\begin{equation}
\label{eq:true_llh}
L(Y|\Omega) =  (2\pi)^{-\frac{np}{2}}\det(\Omega)^{\frac{n}{2}}\exp \left \{-\frac{1}{2} tr(Y \Omega Y^T) \right \}.
\end{equation}

Similar to \cite{atchade2019quasi,zhang2022contraction,sagar2024precision}, we assume the diagonal elements $\omega_{ii}$ of $\Omega$ are known. Without loss of generality, we assume the diagonal $\omega_{ii}=1$ for all $i=1,\ldots,p$. In this work, we consider high-dimensional setting that $ p>n $. We assume that the true precision matrix $ \Omega_0 $ is sparse in the sense that its number of non-zero off-diagonal elements is small.

\subsection{Method}

As in \cite{li2019graphical}, the horseshoe prior \cite{carvalho2010horseshoe} is independently specified on the elements of $ \Omega $ as 
\begin{align}\label{eq:HS}
\omega_{ij} \mid \lambda_{ij},\tau \sim N(0,\lambda_{ij}^2 \tau^2), 
\quad 
\lambda_{ij} \sim \mbox{Ca}_+(0, 1), 
\quad 
\tau \sim \mbox{Ca}_+(0, 1),
\end{align}
for $ 1 \leq i < j \leq p$, where $\mbox{Ca}_+(0, 1)$ denotes the truncated standard half-Cauchy distribution with density proportional to $(1+t^2)^{-1} \mathbbm{1}_{(0, \infty)}(t)$. We shall denote this prior induced by the hierarchy in \eqref{eq:HS} by $\Pi_{HS} $. 

It is important to note that \cite{zhang2022contraction} presents theoretical results exclusively for the case where \( \tau \) is a fixed value. However, in their numerical studies, they use a half-Cauchy prior, as in \cite{li2019graphical}. Therefore, our work extends their theoretical framework.

We focus on the tempered posterior given as
\begin{equation*}
\Pi_{n,\alpha}( \Omega )
\propto
L(Y|\Omega)^\alpha \Pi_{HS} (\Omega)
.
\end{equation*}
for $ \alpha \in (0,1) $. It is also called as fractional posterior \cite{bhattacharya2016bayesian,alquier2020concentration}. The regular posterior is obtained when $ \alpha = 1 $, however, using a smaller value in general would be benifit in both computational \cite{friel2008marginal} or in model mis-specified cases \cite{grunwald2017inconsistency}.  It is worth mentioning that tempered/fractional posteriors has received significant attention in recent years, as in \cite{agapiou2023heavy,bissiri2013general,chakraborty2020bayesian,liu2021bayesian,l2023semiparametric,medina2022robustness,syring2019calibrating,yonekura2023adaptation,yang2020alpha}.

Put $ \hat\Omega = \int \Omega 	\Pi_{n,\alpha}({\rm d} \Omega )  $
 as our mean estimator.

\section{Results}
\label{sc_concentration}
\subsection{Consistency of tempered posterior}

Let $ p_{\Omega} $ and $ P_{\Omega} $ denote the density and the distribution of a Gaussian $ \mathcal{N} (0, \Omega^{-1} ) $ random variable. 
Initially, we enumerate our assumptions regarding the truth.

\begin{assume}[Bounded spectrum]
	\label{assum_true_precision_Spectrum}
	The true precision matrix $\Omega_0 $ satisfies that there	exists a universal constant $ \varepsilon_0 >0 $ that
	$
	 0 < \varepsilon_0^{-1} \leq \eig_{\rm min}(\Omega) 
	\leq \eig_{\rm max}(\Omega) \leq \varepsilon_0 < \infty
	.
	$
\end{assume}

\begin{assume}
	\label{assum_true_precision}
	The true precision matrix $\Omega_0 $ is sparse in the sense that
	$
	 \Omega_0 \in  \mathcal{S}(s) := 
\left\lbrace 
\Omega : \!\! 
\sum_{ 1\leq i < j \leq p}
\mathbf{1}_{(\omega_{ij} \neq 0)} \leq s
\right\rbrace
.
	$
\end{assume}

\begin{assume}[Sparsity]
	\label{assum_sparsity}
It is assumed that $ s<n $.
\end{assume}

Assumption \ref{assum_true_precision_Spectrum} and \ref{assum_true_precision} are standard in both frequentist and Bayesian literature for sparse precision matrix estimation, \cite{cai2016estimating,jankova2017honest,banerjee2015bayesian,sagar2024precision,zhang2022contraction}. Assumption \ref{assum_sparsity} is a sparsity condition often assumed in high-dimensional settings \cite{bellec2018slope,gao2020general,banerjee2021bayesian}.
We now present first a consistency result of the tempered posterior in $ p > n $ setting.

\begin{theorem}
	\label{thm_expectation}
For any \(\alpha \in (0,1)\), assume that Assumptions \ref{assum_true_precision_Spectrum}, \ref{assum_true_precision}, and \ref{assum_sparsity} are satisfied. Given any \(\Omega_0\) such that \(\|\Omega_0\|_\infty \leq C_1\), then
	\begin{equation*}
	\mathbb{E} \left[ \int D_{\alpha}(P_{\Omega},P_{\Omega_0}) 
	\Pi_{n,\alpha}({\rm d} \Omega ) \right]
	\leq 
	C \frac{1+\alpha}{1-\alpha}\varepsilon_n
	,
	\end{equation*}
	where
	$
	\varepsilon_n
	=
 \frac{s\log (p/s)}{n}
	,
	$
	and $ C $ is a universal positive constant depending only on $ \varepsilon_0, C_1 $.
\end{theorem}

Our proof technique leverages the general concentration theory for tempered posteriors as described in \cite{bhattacharya2016bayesian,alquier2020concentration}, which ensures the desired convergence with respect to the $\alpha$-R\'enyi divergence. However, in the context of precision matrices, the Frobenius norm is more interpretable than the $\alpha$-R\'enyi divergence. Using results from  \cite{banerjee2015bayesian} and \cite{van2014renyi}, the results for the Frobenius norm are stated in the following corollary. 

Put $ c_\alpha :=1 $ for $ \alpha \in [0.5,1) $ and  $ c_\alpha := (1-\alpha)/\alpha $ for $ \alpha \in (0, 0.5)  $.

\begin{cor}
	\label{cor_expect_Frobenius}
	Assuming that Theorem \ref{thm_expectation} holds, then we have that
	\begin{equation*}
	\mathbb{E} \left[ \int 
	\| \Omega - \Omega_0 \|_2^2
	\Pi_{n,\alpha}({\rm d} \Omega ) \right]
	\leq 
	C' c_\alpha 
	\frac{1+\alpha}{1-\alpha}\varepsilon_n
	,
	\,
\text{ 	and an application of Jensen's inequality, } 
\,
\mathbb{E} \left[  
\| \hat\Omega - \Omega_0 \|_2^2
 \right]
\leq 
C' c_\alpha 
\frac{1+\alpha}{1-\alpha}\varepsilon_n
,
\end{equation*}
where $ C' >0 $ is a universal constant depending only on $ \varepsilon_0, C_1 $.
\end{cor}

\subsection{Concentration rates}
\label{sc_rslton_distribu}

We now provide our main result concerning the concentration of the tempered posterior relative to the  $\alpha$-R\'enyi divergence of the densities.

  \begin{theorem}
\label{theorem_main}
Given any \(\alpha \in (0,1)\), assume that Assumptions \ref{assum_true_precision_Spectrum}, \ref{assum_true_precision}, and \ref{assum_sparsity} are satisfied. For all \(\Omega_0\) such that \(\|\Omega_0\|_\infty \leq C_1\), then
	\begin{align*}
	\mathbb{P}\left[
	\int D_{\alpha}(P_{\Omega} ,P_{\Omega_0} ) \Pi_{n,\alpha}({\rm d} \Omega )  
	\leq 
	\frac{2(\alpha+1)}{1-\alpha} \varepsilon_n\right] 
	\geq 
	1-\frac{2}{n\varepsilon_n}
	,
	\end{align*}
		where
	$
	\varepsilon_n
	=
	C \frac{s\log (p/s)}{n}
	,
	$
	and $ C $ is a universal positive constant depending only on $ \varepsilon_0, C_1 $.
\end{theorem}

\begin{rmk}
 Our results imply that the concentration rate $ \varepsilon_n $ are adaptive to the unknown sparsity level $ s $ of the true precision matrix. The concentration rate is of order $ s \log(p/s)/n $ is known to be minimax-optimal in high-dimensional setting \cite{bellec2018slope,gao2020general}. Readers may prefer \cite{banerjee2021bayesian} for a recent review of Bayesian methods in different high-dimensional models. It is noted that our results hold with high probability, specifically depending on \( p \) rather than \( n \). More precisely, the probability is greater than or equal to \( 1 - 2(n\varepsilon_n)^{-1} = 1 - 2(C \, s \log (p/s))^{-1} \).
\end{rmk}

The same results apply when using the Frobenius norm that is stated in the following corollary.

\begin{cor}
	\label{cor_concentration_Hellinger}
	As a special case, Theorem \ref{theorem_main} leads to 
	\begin{equation*}
	\mathbb{P}\left[	\int 
	\| \Omega - \Omega_0 \|_2^2
	\, \Pi_{n,\alpha}({\rm d} \Omega )  
	\leq 
2 c_\alpha 
\frac{1+\alpha}{1-\alpha}
	\varepsilon_n \right]
	\geq 
	1-\frac{2}{n\varepsilon_n}
	, 
	\quad
\text{ and } 
\quad
\mathbb{P}\left[ 
\| \hat\Omega - \Omega_0 \|_2^2 
\leq 
2 c_\alpha 
\frac{1+\alpha}{1-\alpha}
\varepsilon_n \right]
\geq 
1-\frac{2}{n\varepsilon_n}
.
\end{equation*}
\end{cor}

Corollary \ref{cor_concentration_Hellinger} shows that the tempered posterior distribution of $ \Omega $ concentrates around its true value, $ \Omega_0 $, at a specified rate relative to the squared Frobenius norm. As follow with previous works on tempered posterior \cite{chakraborty2020bayesian,alquier2020concentration,bhattacharya2016bayesian}, our results do not require that the true value can be tested against sufficiently separated other values in some suitable sieve as in \cite{sagar2024precision} and \cite{zhang2022contraction}.

One may note that our tempered posterior estimates are not necessarily positive definite. In some cases, a positive semi-definite (PSD) matrix may be more preferable. To obtain a PSD matrix \( \tilde{\Omega} \), a post-processing step from \cite{zhang2022contraction} can be used:
\[
\tilde{\Omega} := \argmin_{\Omega \text{ is PSD }} \| \Omega - \bar{\Omega} \|_2,
\]
where \( \bar{\Omega} \) represents random samples from the posterior distribution \( \Pi_{n,\alpha}( \Omega ) \). As shown in \cite{higham1988computing}, this optimization problem has a closed-form solution given by \( \tilde{\Omega} = (B + H)/2 \), where \( H \) is the symmetric polar factor of \( B = (\bar{\Omega} + \bar{\Omega}^\top)/2 \). This leads to the inequality
$
\| \tilde{\Omega} - \Omega_0 \|_2 \leq 2 \| \Omega - \Omega_0 \|_2,
$
which immediately ensures the posterior concentration for the PSD matrix \( \tilde{\Omega} \) from Corollary \ref{cor_concentration_Hellinger}, so the proof is omitted and summarized in the following corollary.

\begin{cor}
Assume that Theorem \ref{theorem_main} holds, then
	\begin{equation*}
	\mathbb{P}\left[	\int 
	\| \tilde\Omega - \Omega_0 \|_2^2
	\, \Pi_{n,\alpha}({\rm d} \Omega )  
	\leq 
2 c_\alpha 
\frac{1+\alpha}{1-\alpha}
	\varepsilon_n \right]
	\geq 
	1-\frac{2}{n\varepsilon_n}
	.
	\end{equation*}
\end{cor}

\subsection{Result in the misspecified case}
\label{sc_mispecifiedcase}
In this section, we show that our previous results can be extended to the misspecified setting.
Assume that the true data generating distribution is parametrized by $ \Omega_0 \notin \mathcal{S}(s) $ and define $P_{\Omega_0} $ as the true distribution. 
Put
\begin{align*}
\Omega_* 
:=
\arg\min_{\Omega \in \mathcal{S}(s)} 
\mathcal{K}(P_{\Omega_0} ,P_{\Omega} )
.
\end{align*}
In order not to change all
the notation, we define an extended parameter set $ \{\Omega_0 \} \cup \mathcal{S}(s) $. Here, we clearly state that only $ \Omega_*  $ is in $ \mathcal{S}(s) $. Example for this circumstance is that $ \Omega_0 $ may no longer be sparse.

\begin{assume}
	\label{assum_pseudotrue_precision}
	The minimal precision matrix $\Omega_*  $ satisfies that
	$$ \Omega_*  \in  \mathcal{S}(s, \varepsilon_0) := 
	\left\lbrace 
	\Omega : \!\! 
	\sum_{ 1\leq i < j \leq p}
	\mathbf{1}_{(\omega_{ij} \neq 0)} \leq s, 
	\, 0 < \varepsilon_0^{-1} \leq \eig_{\rm min}(\Omega) 
	 \leq \eig_{\rm max}(\Omega) \leq \varepsilon_0 < \infty
	\right\rbrace
	.
	$$ 
\end{assume}

\begin{theorem}
	\label{theorem_misspecified}
	For any $\alpha\in(0,1)$, let assume that Assumption  \ref{assum_pseudotrue_precision} and \ref{assum_sparsity} hold. For all \(\Omega_0\) such that \(\|\Omega_*\|_\infty \leq C_1\), then
	\begin{equation*}
	\mathbb{E} \left[ \int 
	D_{\alpha} ( P_{\Omega} ,P_{\Omega_0} ) \Pi_{n,\alpha}({\rm d} \Omega )  \right]
	\leq 
	\frac{\alpha}{1-\alpha} 
\min_{\Omega \in \mathcal{S}(s, \varepsilon_0) } 
\mathcal{K}(P_{\Omega_0} ,P_{\Omega} )
	+ 
	\frac{1+\alpha}{1-\alpha} \varepsilon_n
	,
	\end{equation*}
	where
$
\varepsilon_n
=
\mathcal{C} \frac{s\log (p/s)}{n}
,
$
and $ \mathcal{C} $ is a universal positive constant depending only on $ \varepsilon_0 , C_1 $.
\end{theorem}

In Theorem \ref{theorem_misspecified}, it is important to note that assumptions are made solely regarding \( \Omega_* \) and not \( \Omega_0 \). Consequently, if \( \Omega_0 \in \mathcal{S}(s, \varepsilon_0) \), meaning the model is well-defined and \( \mathcal{K}(P_{\Omega_0} ,P_{\Omega} ) = 0 \) within \( \mathcal{S}(s, \varepsilon_0) \), then Theorem \ref{thm_expectation} can be viewed as a specific instance of Theorem \ref{theorem_misspecified}. Otherwise, the result establishes a general oracle inequality.

Although it is not a sharp oracle inequality due to the differing risk measures on both sides, this observation remains valuable, particularly when $\mathcal{K}(P_{\Omega_0} ,P_{\Omega^*}) $ is minimal. Nonetheless, under additional assumptions, we can further derive an oracle inequality result with $\ell_2$ distance on both sides. The result is as follows.

\begin{cor}
	\label{cor_misspecified}
Assume that Theorem \ref{theorem_misspecified} holds. Then, there exist a universal constant $ \mathcal{C}' >0 $ depending only on $ \varepsilon_0 , C_1 $ such that
	\begin{equation*}
\mathbb{E} \left[ \int 
	\| \Omega - \Omega_0 \|_2^2 \,
 \Pi_{n,\alpha}({\rm d} \Omega )  \right]
\leq 
\mathcal{C}' \left( 
\frac{\alpha}{1-\alpha} 
\min_{\Omega \in \mathcal{S}(s, \varepsilon_0) } 
	\| \Omega - \Omega_0 \|_2^2
+ 
\frac{1+\alpha}{1-\alpha} \varepsilon_n
\right)
,
\end{equation*}
	\begin{equation*}
\mathbb{E} \left[  
\| \hat\Omega - \Omega_0 \|_2^2  \right]
\leq 
\mathcal{C}' \left( 
\frac{\alpha}{1-\alpha} 
\min_{\Omega \in \mathcal{S}(s, \varepsilon_0) } 
	\| \Omega - \Omega_0 \|_2^2
+ 
\frac{1+\alpha}{1-\alpha} \varepsilon_n
\right)
.
\end{equation*}
\end{cor}

To the best of our knowledge, results in Corollary \ref{cor_misspecified} are completely novel for estimating sparse precision matrix from a Bayesian perspective.

\section{Numerical examples}
\label{sc_simu}
In this section, we conduct several simulations to provide a thorough demonstration of our theoretical findings. It is worth noting that more extensive numerical studies involving the graphical horseshoe method, particularly in comparison with other approaches, have been previously presented in the works of \cite{li2019graphical} and \cite{zhang2022contraction}. For the purpose of our analysis, we utilize the Gibbs sampling algorithm, as detailed in \cite{li2019graphical}, to implement our framework. Additionally, the \texttt{R} code for the graphical horseshoe method can be accessed through the \texttt{R} package `\texttt{GHS}' \cite{GHSpackage}. Our \texttt{R} code for simulations is available at \url{https://github.com/tienmt/f_graphical_Horseshoe} .

We consider two distinct settings for the true precision matrix in our analysis:
\begin{itemize}
	\item Sparse setting: In this scenario, the dimension of the precision matrix is set at \( p = 100 \), with a sample size of \( n = 30 \). The true precision matrix \( \Omega_0 \) is characterized by sparsity, having \( s = 86 \) non-zero off-diagonal elements.
	\item Dense setting: Here, the dimension of the precision matrix is \( p = 50 \), with a sample size of \( n = 30 \). The true precision matrix \( \Omega_0 \) is defined as the inverse of the covariance matrix \( \Sigma \), where the diagonal entries of \( \Sigma \) are all equal to 1, and each off-diagonal entry is 0.2, reflecting a denser structure.
\end{itemize}

In addition to considering the Gaussian graphical model, where the data \( Y_i \sim \mathcal{N}(0, \Omega_0^{-1}) \), we also explore a scenario in which the data is generated from a multivariate Student-t distribution with 3 degrees of freedom, denoted as \( Y_i \sim {\rm mvT}_3(0, \Omega_0^{-1}) \). For all settings, the Gibbs sampler is executed for a total of 1100 iterations, with the initial 100 iterations discarded as burn-in steps. Each simulation setting is independently repeated 50 times to ensure robust results. 

We consider three distinct values for the tempered/fractional parameter \( \alpha \) in our analysis: 0.9, 0.5, and 0.1. We refer to our proposed method corresponding to these values as `fGHS$_{\alpha = 0.9}$', `fGHS$_{\alpha = 0.5}$', and `fGHS$_{\alpha = 0.1}$', respectively. For comparison, we denote the method from \cite{li2019graphical}, which corresponds to \( \alpha = 1 \), as `GHS'. The average squared error for estimating the precision matrix, along with the corresponding standard deviations, are reported in Table \ref{tb_1_ttm}.


\input{table_STAT-24-0212_R1.tex}

The results provided in Table \ref{tb_1_ttm} offer a comprehensive evaluation of our proposed method, which leverages the tempered posterior to handle model misspecification effectively. These findings suggest that our method remains robust even when the underlying model is not perfectly specified. Specifically, in situations where the data are generated from the true Gaussian model, our approach exhibits behavior that closely mirrors that of the conventional Bayesian method, particularly when the tempered/fractional parameter \( \alpha \) is near 1. This similarity underscores the consistency of our method with traditional approaches under ideal conditions.

Moreover, our simulations highlight that certain intermediate values of \( \alpha \), such as \( \alpha = 0.9 \) or \( \alpha = 0.5 \), can potentially enhance performance, yielding better results than when \( \alpha = 1 \) is used. This observation points to the potential benefits of fine-tuning the parameter \( \alpha \) to optimize outcomes in practical applications, see \cite{wu2023comparison} for more discussions. However, it is also evident from our studies that selecting a small value of \( \alpha \), e.g. 0.1, can lead to a noticeable increase in error, indicating that there is a delicate balance to be maintained when choosing this parameter.

Overall, these simulation studies provide clear and convincing evidence in support of our theoretical framework. They not only validate the robustness of our method but also demonstrate its flexibility in adapting to different scenarios, thereby offering valuable insights into its practical utility.

\section{Conclusion}
\label{sc_conclusion}
This study has investigated the concentration properties of the tempered posterior when estimating sparse precision matrices using the horseshoe prior. Unlike prior studies that employed either a horseshoe-like prior or a simplified horseshoe prior with a fixed shrinkage parameter, our analysis focuses on the fully specified horseshoe prior. We extend the previous theoretical works by considering a high-dimensional setting ($ p > n $). 
We also present a novel aspect by deriving theoretical results for model misspecification, establishing a general oracle inequality for the posterior in an abstract setting.

We note that our findings are primarily derived using the Frobenius norm, so extending these results to other norms represents a valuable avenue for future research. Another potential gap in our work is that establishing a sharp oracle inequality for the model misspecified case remains an open question for future research.

\bmsection*{Acknowledgments}
The author wishes to acknowledge the support provided by the Centre for Geophysical Forecasting, funded by the Norwegian Research Council under grant no. 309960, at NTNU. Special thanks are extended to the associate editor and the anonymous reviewer for their insightful feedback, which greatly improved the manuscript.

\bmsection*{Conflict of interest}

No, there is no conflict of interest.

\bmsection*{Data Availability Statement}

Data sharing is not applicable to this article as no datasets were generated or analyzed during the current study.

\appendix
\section{Proofs}
\label{sc_proofs}

\subsection{Proof of Theorem~\ref{thm_expectation}}

\begin{proof}	
	We can check the hypotheses on the KL between the likelihood terms as required in Theorem 2.6 in \cite{alquier2020concentration}.
	For 
	$$ \mathcal{K} (p_1,p_2) = \int p_1\log(p_1/p_2)
	,
	$$  
	and let $d_1,\ldots,d_p$ denote the eigenvalues of $\Omega_0^{-1/2}\Omega\Omega^{-1/2}_0$. Then, using Lemma~\ref{lemma:KL}, we have,
	\begin{eqnarray}
	\mathcal{K} (p_{\Omega_0},p_{\Omega}) 
	=
	-\dfrac{1}{2}\sum_{i=1}^{p}\log d_i - \dfrac{1}{2}\sum_{i=1}^{n}(1 - d_i),\nonumber
	.
	\end{eqnarray}
	As $\sum_{i=1}^{n}(1 - d_i)^2  = \|I_p - \Omega_0^{-1/2}\Omega\Omega^{-1/2}_0\|_2^2,$  when  $\|I_p - \Omega_0^{-1/2}\Omega\Omega^{-1/2}_0\|_2^2$ is small, we have, $\max_{1 \leq i \leq p}|1-d_i| < 1$, see \cite{banerjee2015bayesian}. This leads to that,  
	\begin{align*}
	\mathcal{K} (p_{\Omega_0},p_{\Omega})
	\lesssim 
	\sum_{i=1}^n(1-d_i)^2
	= 
	\|I_p - \Omega_0^{-1/2}\Omega\Omega^{-1/2}_0\|_2^2 
	= 
	\|\Omega_0^{-1/2}( \Omega-\Omega_0) \Omega^{-1/2}_0\|_2^2  
	\leq 
	\|\Omega_0^{-1}\|_2^2\|\Omega - \Omega_0\|_2^2 
	\leq 
	\varepsilon_0^{-2} \|\Omega - \Omega_0\|_2^2
	.
	\end{align*}
	Integrating with respect to $\rho_n \propto \mathbf{1}_{\|\Omega - \Omega_0\|_2 < \delta } \Pi_{HS} $ where $\delta= [s\log (p/s)/n ]^{1/2} $, we have that
	\begin{align}
	\label{eq_boundKL}
	\int \mathcal{K} (p_{\Omega_0},p_{\Omega}) 
	\rho_n({\rm d} \Omega)
	\leq 
	2\varepsilon_0^{-2} \delta^2
	\leq 
	2\varepsilon_0^{-2} 
	\frac{s\log (p/s)}{n}
	.
	\end{align}	
	Applying Lemma \ref{lm_bound_prior_horseshoe} with $ d = p(p-1)/2 $, we have that 
	\begin{align*}
	\frac1n \mathcal{K}(\rho_n,\pi)
	\leq
	\frac1n \log \frac{1}{	\Pi_{HS} (\|\Omega - \Omega_0 \|_2<\delta)}
	\leq
	\frac{1}{n} K s\log (p/s)
	.
	\end{align*}
	Consequently, we can now apply Theorem 2.6 from \cite{alquier2020concentration} to obtain the result, with
$
	\varepsilon_n
	=
	2\varepsilon_0^{-2} \frac{s\log (p/s)}{n}
	.
$
	The proof is completed.	
\end{proof}

\begin{proof}[\bf Proof of Corollary \ref{cor_expect_Frobenius}]
	From \cite{van2014renyi}, we have that 
	$$ 
	h^2(p_{\Omega},p_{\Omega_0})
	\leq 
	D_{1/2}(P_{\Omega},P_{\Omega_0})  
	\leq 
	D_{\alpha}(P_{\Omega},P_{\Omega_0}) 
	,
	$$
	for $ \alpha \in [0.5,1) $. In addition, from \cite{van2014renyi}, all $\alpha$-R\'enyi divergences are all equivalent for $0<\alpha<1$, through the formula $ \frac{\alpha}{\beta}\frac{1-\beta}{1-\alpha} D_{\beta} \leq D_\alpha \leq D_\beta $ for $\alpha \leq \beta $, thus we have that 
	$$ 
	D_{1/2}(P_{\theta},P_{\theta_0})  
	\leq 
	\frac{(1-\alpha)1/2}{\alpha (1-1/2)} D_{\alpha}(P_{\theta},P_{\theta_0}) 
	=
	\frac{(1-\alpha)}{\alpha} D_{\alpha}(P_{\theta},P_{\theta_0}) 
	,
	$$ 
	for $ \alpha \in (0, 0.5)  $.
	Thus, with $ c_\alpha :=1 $ for $ \alpha \in [0.5,1) $ and  $ c_\alpha := (1-\alpha)/\alpha $ for $ \alpha \in (0, 0.5)  $, we obtain from Theorem \ref{thm_expectation} that 
	\begin{equation*}
	\mathbb{E} \left[ \int 	
	\frac{1}{c_\alpha}
	h^2 (p_{\Omega},p_{\Omega_0}) \pi_{n,\alpha}({\rm d} \Omega ) \right]
	\leq 
	\frac{1+\alpha}{1-\alpha}\varepsilon_n
	.
	\end{equation*} 
	From Lemma \ref{lm_Froben_hellinger}, one gets that
	\begin{equation*}
	\mathbb{E} \left[ \int 	
	\frac{c_0}{c_\alpha}
	\| \Omega_1 - \Omega_2 \|_2^2
	\pi_{n,\alpha}({\rm d} \Omega ) \right]
	\leq 
	\frac{1+\alpha}{1-\alpha}\varepsilon_n
	.
	\end{equation*}
	The proof is completed.
\end{proof}

\subsection{Proof of Theorem~\ref{theorem_main}}

\begin{proof}
	\label{proof_theorem_main}		

We utilize the general framework of posterior concentration rates by confirming the criteria outlined in Theorem 2.4 of \cite{alquier2020concentration}. This involves assessing the prior concentration rate within Kullback–Leibler neighborhoods.
	
	Put 
	$$ 
	\mathcal{K} (p_1,p_2) = \int p_1\log(p_1/p_2), 
	\text{ and }
 V(p_1,p_2) = \int p_1 \log^2(p_1/p_2)
 .
 $$ and let $ d_1,\ldots,d_p $ denote the eigenvalues of $ \Omega_0^{-1/2}\Omega\Omega^{-1/2}_0 $. Controlling the term $ \mathcal{K} (p_1,p_2) = \int p_1\log(p_1/p_2) $ is presented in the proof of Theorem \ref{thm_expectation}, given in page \pageref{eq_boundKL}.	Now, using Lemma~\ref{lemma:KL}, we have,
	\begin{eqnarray*}
	V(p_{\Omega_0},p_{\Omega}) 
	= 
	\dfrac{1}{2}\sum_{i=1}^{n}(1 - d_i)^2 + \mathcal{K} (p_{\Omega_0},p_{\Omega})^2.
	\end{eqnarray*}
	As $\sum_{i=1}^{n}(1 - d_i)^2  = \|I_p - \Omega_0^{-1/2}\Omega\Omega^{-1/2}_0\|_2^2,$  when  $\|I_p - \Omega_0^{-1/2}\Omega\Omega^{-1/2}_0\|_2^2$ is small, we have, $\max_{1 \leq i \leq p}|1-d_i| < 1.$ This gives 
	\begin{align*}
	V(p_{\Omega_0},p_{\Omega}) 
	\lesssim 
	\sum_{i=1}^n(1-d_i)^2
	= 
	\|I_p - \Omega_0^{-1/2}\Omega\Omega^{-1/2}_0\|_2^2 
	= 
	\|\Omega_0^{-1/2}( \Omega-\Omega_0) \Omega^{-1/2}_0\|_2^2  
	\leq  
	\|\Omega_0^{-1}\|_2^2\|\Omega - \Omega_0\|_2^2 
	\leq  
	\varepsilon_0^{-2}\|\Omega - \Omega_0\|_2^2
	.
	\end{align*}
	
	Integrating with respect to $\rho_n \propto \mathbf{1}_{\|\Omega - \Omega_0\|_2 < \delta } \Pi_{HS} $ where $\delta= [s\log (p/s)/n ]^{1/2} $, we have that
	\begin{align*}
	\int	V(p_{\Omega_0},p_{\Omega}) 
	\rho_n({\rm d} \Omega)
	\leq 
	2\varepsilon_0^{-2} \delta^2
	\leq 
	2\varepsilon_0^{-2} 
	\frac{s\log (p/s)}{n}
	.
	\end{align*}
	and from \eqref{eq_boundKL},
$
	\int	\mathcal{K} (p_{\Omega_0},p_{\Omega}) 
	\rho_n({\rm d} \Omega)
	\leq 
	2\varepsilon_0^{-2} \delta^2
	\leq 
	2\varepsilon_0^{-2} 
	\frac{s\log (p/s)}{n}
	$.		
	From Lemma \ref{lm_bound_prior_horseshoe}, we have that 
	\begin{align*}
	\frac1n \mathcal{K}(\rho_n,\pi)
	\leq
	\frac1n \log \frac{1}{	\Pi_{HS} (\|\Omega  - \Omega_0\|_2<\delta)}
	\leq
	\frac{1}{n} K s\log (p/s)
	.
	\end{align*}
	Consequently, we can now apply Theorem 2.4 from \cite{alquier2020concentration} to obtain the result,  with
$
	\varepsilon_n
	=
	2\varepsilon_0^{-2} 
	\frac{s\log (p/s)}{n}
	.
$
	The proof is completed.
\end{proof}

\begin{proof}[\bf Proof of Corollary~\ref{cor_concentration_Hellinger}]
	The proof is similar to the proof of Corollary \ref{cor_expect_Frobenius}.
\end{proof}

\begin{lemma}[Lemma B1 in \cite{sagar2024precision}]
	\label{lemma:KL}
	Let $p_k$ denote the density of a $\mathcal{N}_d(0, \Sigma_k)$ random variable, $k = 1,2.$ Denote the corresponding precision matrices by $\Omega_k = \Sigma_k^{-1},k=1,2.$ Then,
	\begin{align*}
		\mathbb{E}_{p_1}\left\lbrace \log \dfrac{p_1}{p_2}(X)\right \rbrace  
		&= 
		\dfrac{1}{2}\left\lbrace \log \det \Omega_1 - \log \det \Omega_2 + trace(\Omega^{-1}_1\Omega_2 - I_d) \right\rbrace,  
		\\
		{\rm Var}_{p_1}\left\lbrace \log \dfrac{p_1}{p_2}(X)\right \rbrace 
		&=
		\dfrac{1}{2}\, trace\{(\Omega_1^{-1/2}\Omega_2\Omega_1^{-1/2} - I_d)^2\}. 
	\end{align*}
\end{lemma}

\begin{lemma}[Lemma A.1 in \cite{banerjee2015bayesian}]
	\label{lm_Froben_hellinger}
	Let $p_k$ denote the density of a $\mathcal{N}_d(0, \Sigma_k)$ random variable, $k = 1,2.$ Denote the corresponding precision matrices by $\Omega_k = \Sigma_k^{-1},k=1,2.$ Then, there exists a positive constant $ c_0 $ such that
$
\| \Omega_1 - \Omega_2 \|_2^2
	\leq
c_0 h^2 ( p_1 , p_2 )
.
$
\end{lemma}

The subsequent lemma introduces a novel prior concentration result for the horseshoe prior, offering a lower bound on the probability assigned to a Euclidean neighborhood of a sparse vector. This lemma represents a vital improvement over Lemma 1 in \cite{chakraborty2020bayesian}, crucial for achieving the optimal rate \( s\log(p/s)/n \) instead of the sub-optimal rate \( s\log(p)/n \) in high dimensional setting, see \cite{bellec2018slope}, \cite{gao2020general}.
\begin{lemma}
	\label{lm_bound_prior_horseshoe}
	Suppose $ \theta_0 \in \mathbb{R}^d $ such that $ \|\theta_0\|_0 = s $ and  that $ s < n < d $ and $ \|\theta_0\|_\infty \leq C_1  $. Suppose $\theta \sim \Pi_{HS} $. Define $\delta=\{s\log (d/s)/n\}^{1/2}$. Then we have, for some constant $K>0 $, that
	$
	\Pi_{HS} (\|\theta -\theta_0\|_2<\delta)
	\geq 
	e^{-Ks\log (d/s)}.
	$
\end{lemma}

\begin{proof}[Proof of Lemma \ref{lm_bound_prior_horseshoe}]
	Using the scale-mixture formulation of $\Pi_{HS}$,
	\begin{equation*}
	\begin{split}
	\Pi_{HS} (\| \theta -\theta_0\|_2<\delta)
	=
	\int_\tau pr(\| \theta -\theta_0\|_2<\delta\mid \tau) f(\tau) {\rm d}\tau
 \geq 
	\int_{I_{\tau_*}} pr(\| \theta -\theta_0\|_2<\delta\mid \tau) f(\tau) {\rm d}\tau,
	\end{split}
	\end{equation*}
	where $I_{\tau_*}=[\tau_*/2,\tau_*]$ with $\tau_* = s(s/p)^{3/2}\{\log( d/s)/n\}^{1/2}$. Let $S=\{1\leq j \leq d : \theta_{0j}\neq 0\}$.
	For $\tau \in I_{\tau^*}$, one can lower bound the conditional probability as
	\begin{equation}\label{eq:split_prob}
	\begin{split}
	pr(\| \theta -\theta_0\|_2 <\delta\mid \tau)
	&\geq 
	pr\left(\|\theta_S-\theta_{0S}\|_2<\delta/2 \mid \tau\right)
	pr\left(\| \theta_{S^c}\|_2<\delta/2 \mid \tau\right)
	\\
	&\geq \prod_{j\in S}pr\left(|\theta_j-\theta_{0j}|
	<
	\frac{\delta}{2s}\mid \tau\right) 
	\prod_{j \in S^c}pr\left(|\theta_j|<\frac{\delta}{2\sqrt{d}}\mid \tau\right).
	\end{split}
	\end{equation}
	For a fixed $\tau \in I_{\tau_*}$, we will  lower bounds each of the terms in the right hand side of \eqref{eq:split_prob}. 
	
	Initially, we examine \( pr\{ \, | \theta_j | <\delta /(2d^{1/2})\mid \tau\} \) with \( \tau \) in \( I_{\tau_*} \). Given \( \tau \) and \( \lambda \), \( \theta_j \) follows a normal distribution \( N (0,\lambda_j^2\tau^2) \). Leveraging the Chernoff bound for a Gaussian random variable, we obtain:
	\begin{align*}
	pr\big\{ |\theta_j| > \delta/(2d^{1/2}) \mid \lambda_j, \tau\big \}
	\leq 
	2e^{-\delta^2/(8p\lambda_j^2 \tau^2)}
	\leq 
	2e^{-\delta^2/(8p\lambda_j^2 \tau_*^2)}
	.
	\end{align*}
	Thus, 
	\begin{align*}
	pr\big\{ | \theta_j | <\delta/(2d^{1/2}) \mid \tau \big \}
	 =
	\int_{\lambda_j} pr \big\{ |\theta_j| < \delta/(2d^{1/2}) \mid \lambda_j, \tau\big\} \, f(\lambda_j) \,{\rm d}\lambda_j
	& \geq \int_{\lambda_j} \left\lbrace1-2\exp{\left( -
		\delta^2/(8d\lambda_j^2 \tau_*^2) 
		\right)}\right\rbrace  \, f(\lambda_j) \,{\rm d}\lambda_j
	\\
	&= 1-\dfrac{4}{\pi}\int_{\lambda_j} 
	\exp{\left(
		-\delta^2/(8d\lambda_j^2 \tau_*^2)
		\right)} (1+\lambda_j^2)^{-1}{\rm d}\lambda_j=1-\dfrac{4}{\pi} \, \mathcal{I}
	.
	\end{align*}
	We now upper-bound the integrand as follows,
	\begin{align*}
	\mathcal{I}
&		=
	\int_{\lambda_j} \exp{\left(-\dfrac{\delta^2}{8d\tau_*^2 \lambda_j^2} \right)} \, (1+\lambda_j^2)^{-1}{\rm d}\lambda_j 
\\
&	\leq 
	\int_{\lambda_j} \exp{\left(-\dfrac{\delta^2}{8d\tau_*^2\lambda_j^2}\right)}\lambda_j^{-2} {\rm d}\lambda_j
 = 
	\frac{1}{2} \int_0^{\infty} z^{-1/2}  \exp{\left(-\dfrac{\delta^2z}{8d\tau_*^2}\right)}
	{\rm d}z
 = 
	\dfrac{\Gamma (1/2)}{\{2 \delta^2/(8d \tau_*^2)\}^{1/2}} 
	=  
	\dfrac{ 2s^2\sqrt{ \pi } }{d} 
	,
	\end{align*}
	where a substitution $z=1/\lambda^2$ at the third step was made. Thus, for $\tau \in I_{\tau_*}$, 
	$
	pr(| \theta_j | < \delta/2d^{1/2}\mid \tau)
	\geq 
	1 - \dfrac{8s^2}{d\sqrt{\pi}}
	.
	$
	\\	
	Now, we lower bound $pr(| \theta_j-\theta_{0j} | <\delta_0|\tau)$ with $\tau \in I_{\tau_*}$. Letting $\delta_0 = \delta/(2\sqrt{s}) $, then
	\begin{align*}
	\begin{split}
	pr( | \theta_j-\theta_{0j}| <\delta_0 \mid \tau)
	& =
	(\frac{2}{\pi^3})^{1/2}\int_{\lambda_j}\int_{\mid \theta_j-\theta_0 \mid <\delta_0} \exp\{{-\theta_j^2/(2\lambda_j^2\tau^2)}\} \dfrac{1}{\lambda_j\tau(1+\lambda_j^2)} \,{\rm d}\theta_j {\rm d}\lambda_j
	\\
	& \geq (\frac{2}{\pi^3})^{1/2}\int_{|\theta_j-\theta_0|<\delta_0} \int_{1/\tau}^{2/\tau} \exp\{{-\theta_j^2/(2\lambda_j^2\tau^2)}\} \dfrac{1}{\lambda_j\tau(1+\lambda_j^2)} \,{\rm d}\lambda_j {\rm d}\theta_j
	\\
	& \geq 
	(\frac{2}{\pi^3})^{1/2} \int_{|\theta_j-\theta_0|<\delta_0} \exp(-\theta_j^2/2)  \left(\int_{1/\tau}^{2/\tau} \dfrac{1}{1+\lambda_j^2} \,{\rm d}\lambda_j \right) {\rm d}\theta_j,
	\\
	\end{split}
	\end{align*}
	since for $\lambda_j \in [1/\tau, 2/\tau]$, $1/(\lambda_j \tau) \geq 1/2$ and $\exp\{{-\theta_j^2/(2\lambda_j^2\tau^2)}\} \geq \exp{(-\theta_j^2/2)}$. Continuing,
	\begin{align*}
	pr( | \theta_j-\theta_{0j} | <\delta_0\mid \tau)
&	\geq 
	(2/\pi^3)^{1/2}\, \dfrac{\tau}{4+\tau^2}\, \int_{|\theta_j-\theta_0|<\delta_0} \exp{(-\theta_j^2/2)}{\rm d}\theta_j
\\	& \geq 
	(2/\pi^3)^{1/2}\,
	\dfrac{\tau}{4+\tau^2}\,
	 \exp\{{-( C_1 + 1)^2/2}\}\,\delta_0
	\\
	&\geq 
	K \, \tau \, \delta_0
	\geq
	K \, \frac{\tau_* \delta}{4 \sqrt{s} }
	\geq
	K	s \left(\frac{s}{p}\right)^{3/2} \frac{\log (p/s)}{n}
	\geq
	K	 \left(\frac{s}{p}\right)^{5/2}
	=
	K	e^{-(5/2)\log(p/s)}
	,
	\end{align*}
	where in the third step, we used $4+\tau^2<5$ and $ n<p $ in the final step.
	
	By substituting these bounds into \eqref{eq:split_prob}, we obtain, for \( \tau \in I_{\tau_*} \),
	\begin{align*}
	pr(\|\theta-\theta_0\|_2<\delta \mid \tau )
	\geq
	\left( K	e^{-(5/2)\log(p/s)} \right)^s 
	\left(	1 - \dfrac{8s^2}{p\sqrt{\pi}}  \right)^{p-s}
	\geq
	e^{-Ks\log (p/s)},
	\end{align*}
	where $K$ is a positive constant. The proof concludes by noting that the probability of $ \tau $ being within the interval \(I_{\tau_*}\) is greater than or equal to \( \tau_*/(2\pi)\). Thus, with a minor deviation in notation, we obtain,
	\begin{equation*}
	\Pi_{HS}(\|\theta-\theta_0\|_2<\delta) \geq e^{-Ks\log (p/s)},
	\end{equation*} 
	for some positive constant $K$. The proof is completed.
\end{proof}

\subsection{Proof of Theorem~\ref{theorem_misspecified}}

\begin{proof}
	
	We can check the assumptions as in Theorem 2.7 in \cite{alquier2020concentration}.
	
	Let us define $ A = \Omega_1^{-1/2}\Omega_2\Omega_1^{-1/2}.$ Note that, for a random variable $ Z \sim \mathcal{N}_d(0,\Sigma)$, we have,
	$\mathbb{E}( Z^T A Z) = trace(A \Sigma)
	.$
	Then, for $ X \sim \mathcal{N}_d(0,\Sigma_0)$,
	\begin{align*}
	\mathbb{E}_{p_0}\left\lbrace \log \dfrac{p_1}{p_2}(X)\right \rbrace  
	=
	\dfrac{1}{2} \left\lbrace \log \det \Omega_1 - \log \det \Omega_2 + \mathbb{E}_{p_0}(X^T(\Omega_2 - \Omega_1)X)\right\rbrace  
=
	\dfrac{1}{2} \left\lbrace \log \det \Omega_1 - \log \det \Omega_2 + trace[(\Omega_2 - \Omega_1)\Sigma_0] \right\rbrace 
	.
	\end{align*}
	Thus, we have that
	\begin{align*}
	\mathbb{E}_{\Omega_0} \left[\log
	\frac{ p_{\Omega_*}}{ p_{\Omega} } (X) \right]
		\leq
	\dfrac{1}{2} \left\lbrace \log \left( \frac{\det (\Omega_*)}{\det (\Omega) } 
	\right)
	+ 
	trace[(\Omega - \Omega_*)\Omega_0^{-1}] \right\rbrace 
&	\leq
	\dfrac{1}{2} \left\lbrace \log \left( \frac{\det (\Omega_*)}{\det (\Omega) } 
	\right)
	+ 
	\| (\Omega - \Omega_*)\Omega_*^{-1} \Omega_*\|_2 \| \Omega_0^{-1}\|_2 \right\rbrace 
	\\
	&	\leq
	\dfrac{1}{2} \left\lbrace 
	\log \left( \frac{\det (\Omega_*)}{\det (\Omega) } 
	\right)
	+ 
	\| \Omega \Omega_*^{-1} - I_p \|_2 \| \Omega_*\|_2  \varepsilon_0^{-1} \right\rbrace 
	.	
	\end{align*}
	Let $ d_1, \ldots, d_p $ denote the eigenvalues of $  \Omega \Omega_*^{-1} $. As $\sum_{i=1}^{n}(1 - d_i)^2  = \|I_p -  \Omega \Omega_*^{-1} \|_2^2,$  when  $\|I_p -  \Omega \Omega_*^{-1}\|_2^2$ is small, we have, $\max_{1 \leq i \leq p}|1-d_i| < 1.$ This leads to that,  
	\begin{align*}
	\mathbb{E}_{\Omega_0} \left[\log
	\frac{ p_{\Omega_*}}{ p_{\Omega} } (X) \right]
		\leq
	\dfrac{1}{2} \left\lbrace 
	\sum_{i=1}^{p} \log d_i
	+ 
	\| \Omega_*\|_2  \varepsilon_0^{-1}
	\sum_{i=1}^{p} ( 1- d_i)
	\right\rbrace 
	\leq  
	K\| \Omega_*\|_2  \varepsilon_0^{-1}
	\sum_{i=1}^{p} ( 1- d_i)^2
	\leq 
	K_{\varepsilon_0}
	\| \Omega - \Omega_* \|_2^2
	,
	\end{align*}
	for some absolute constant $ K_{\varepsilon_0}>0 $ depending only on $ \varepsilon_0 $. 	Integrating with respect to $\rho_n \propto \mathbf{1}_{\|\Omega - \Omega_* \|_2 < \delta } \Pi_{HS} $ where $\delta= [s\log (p/s)/n ]^{1/2} $, we have that
	\begin{align*}
	\int		\mathbb{E}_{\Omega_0} \left[\log
	\frac{ p_{\Omega_*}}{ p_{\Omega} } (X) \right]
	\rho_n({\rm d} \Omega)
	\leq 
	K_{\varepsilon_0} \delta^2
	\leq 
	K_{\varepsilon_0}
	\frac{s\log (p/s)}{n}
	.
	\end{align*}

	To apply Theorem 2.7 in \cite{alquier2020concentration} it remains to compute the KL between the approximation $ \rho_n $ and the prior. 	From Lemma \ref{lm_bound_prior_horseshoe}, we have that 
	\begin{align*}
	\frac1n \mathcal{K}(\rho_n,\pi)
	\leq
	\frac1n \log \frac{1}{	\Pi_{HS} (\|\Omega - \Omega_0\|_2<\delta)}
	\leq
	\frac{1}{n} K s\log (p/s)
	.
	\end{align*}
	To obtain an estimate of the rate $\varepsilon_n$ as in Theorem 2.7 in \cite{alquier2020concentration}, we put together those bounds and choosing 
$
	\varepsilon_n
	=
	K_{\varepsilon_0} \frac{s\log (p/s)}{n}
	.
$
	The results are followed and the proof is completed.	
\end{proof}

\begin{proof}[\bf Proof of Corollary~\ref{cor_misspecified}]
	A similar argument as in the proof of Corollary \ref{cor_expect_Frobenius}, we have that 
	$ 	
	h^2 (p_{\Omega},p_{\Omega_0}) 
	\leq c_\alpha 	D_{\alpha} ( P_{\Omega} ,P_{\Omega_0} )
	$. Then, using Lemma \ref{lm_Froben_hellinger}, one has that $ \| \Omega_1 - \Omega_2 \|_2^2
	\lesssim
 h^2 ( p_1 , p_2 ) $. Additionally, from the proof of Theorem \ref{theorem_main}, page \pageref{eq_boundKL}, one has that
	$  	
	\mathcal{K} (p_{\Omega_0},p_{\Omega})
	\lesssim 	
	\varepsilon_0^{-2} \|\Omega - \Omega_0\|_2^2 
	$. Putting all together into Theorem \ref{theorem_misspecified}, we get that
	\begin{equation*}
\mathbb{E} \left[ \int 
\| \Omega - \Omega_0 \|_2^2 \,
\Pi_{n,\alpha}({\rm d} \Omega )  \right]
\leq 
C \left(
\frac{\alpha}{1-\alpha} 
\min_{\Omega \in \mathcal{S}(s, \varepsilon_0) } 
	\| \Omega - \Omega_0 \|_2^2
+ 
\frac{1+\alpha}{1-\alpha} \varepsilon_n
\right)
,
\end{equation*}
for some universal constant $ C>0 $ that do not depend on $ n, s, p $. Subsequently, an application of Jensen's inequality yields the result for the mean estimator. The proof is completed.
\end{proof}

\end{document}

%% file: table_STAT-24-0212_R1.tex
\begin{table}
	\centering
	\begin{tabular}{ l | c | c | c | c } 
		\hline
Method		& \multicolumn{2}{c|}{dense, $ p=50 , n=30 $}
		& \multicolumn{2}{c}{sparse, $ p=100 , n=30 $}
		\\
		\hline 
 & $ Y_i \sim \mathcal{N} (0, \Omega^{-1}) $ & $ Y_i \sim {\rm mvT}_3 (0, \Omega^{-1}) $ 
& $ Y_i \sim \mathcal{N} (0, \Omega^{-1}) $ & $ Y_i \sim {\rm mvT}_3 (0, \Omega^{-1}) $ 
 \\
		\hline
GHS &  0.013 (0.003) & 0.009 (0.004) 
	& 0.0028 (0.0005) & 0.0041 (0.0019)
		\\ 
fGHS$_{\alpha = 0.9} $ & 0.011 (0.003) & 0.007 (0.003)
	& 0.0024 (0.0004) & 0.0032 (0.0011)
		\\ 
fGHS$_{\alpha = 0.5} $ & 0.003 (0.001) & 0.005 (0.001)
	& 0.0025 (0.0003) & 0.0031 (0.0006)
		\\ 
fGHS$_{\alpha = 0.1} $ & 0.015 (0.005) & 0.006 (0.002)
	& 0.0088 (0.0009) & 0.0042 (0.0006)
		\\
		\hline
	\end{tabular}
	\caption{Estimation errors are measured using the mean squared error, calculated across 50 simulated repetitions.}
	\label{tb_1_ttm}
\end{table}